\documentclass[11pt]{article}
\usepackage{amsmath,amsthm,amsfonts,amssymb,amscd,latexsym}
\usepackage[latin1]{inputenc}
\usepackage[all]{xy}
\usepackage{xypic}
\usepackage{url}
\usepackage{xcolor}
\usepackage{tikz}
\usepackage{tikz-qtree}
\usetikzlibrary{trees}
\usetikzlibrary[trees]
\newtheorem{theorem}{Theorem}
\theoremstyle{remark}

\newtheorem{remark}{Remark}

\newtheorem{lemma}{Lemma}[section]

\theoremstyle{plain}

\title{Reflections on Russell's antinomy}
\author{Paola Cattabriga \\ \small{University of Bologna}}

\date{} 
\begin{document}
\maketitle
\begin{abstract}
We present Russell's antinomy using three distinct deductive systems, which are then compared to deepen the logical deductions that lead to the contradiction.  Some inferential paths are then presented, alternative to the commonly accepted one, that allow for the formal assertion of the antinomy without deriving the contradiction, thus preserving the coherence of the system.  In light of this, the purpose of this article is to propose a review of the consequences of asserting Russell's antinomy and, by extension, the widespread belief that any attempt to resolve a paradox is doomed to failure. 
\end{abstract}

\section{Introduction}
It is not always impossible to demonstrate that a proposition that leads to a contradiction within a system, can also be used to demonstrate that the contradiction does not actually occur within the framework of the system's rules of inference. Said another way, this does not imply that the same system cannot provide alternative inferential paths that preserve coherence.
In this paper, we first introduce three distinct approaches to demonstrating how Russell's antinomy results in contradiction: the traditional argument from Zermelo-Frankel set theory, a first-order logic theorem, and the use of semantic trees. 
Comparing them, the following sections explore other possible inferential paths, their underlying involvement with Leibniz's indiscernibility, and finally a review in light of the Theory of Definition.
The first one we recall is the most simple, in Zermelo-Fraenkel set theory (\cite{foundations}, 31).
\begin{theorem}\label{RA}
There exists no set  which 
contains exactly those elements which do not contain 
themselves, in symbols  $\neg\exists y\forall x(x\in y \iff  x\notin 
x) $.
\end{theorem}
\begin{proof}
By contradiction. Assume that $y$ is a set such that for every 
element 
$x$, $ x \in y$ if and only if $ x \notin x$. For $x=y$, we have $y \in y$
 if and only 
if $y \notin y$. Since, obviously, $y \in y$ or $y \notin y$, and as we 
saw, 
each of $y \in y$ and $y \notin y$ implies the other statement, we have 
both $y 
\in y$ and $y \notin y$, which is a contradiction.
\end{proof}
Next, let us present a development of the same argument using first order logic, which provides deeper insight into the internal logical interconnections of the antinomy\footnote{For supplying the following version of the antinomy, many thanks to Giacomo Lenzi.}.
\begin{theorem}\label{gl}
The following first order formula is provably false:
\begin{equation}\label{eq:1} 
\exists r \forall x ( x\in r\iff x\notin x).
\end{equation}
\end{theorem}
\begin{proof} 
Suppose that (\ref{eq:1}) holds, we derive a contradiction.
In first-order logic, the axiom (or scheme of axioms, see the Predicative Axiom in \cite{CK}) holds
\begin{equation}\label{eq:2} 
\forall x \phi(x,r) \to\phi(r,r)
\end{equation}
provided that the replacement of  $x$ with $r$  is sound ( but it certainly is if $\phi$ has no quantifiers).
By the generalization rule it then follows
\begin{equation}\label{eq:3}
 \forall r(\forall x\phi(x,r) \to\phi(r,r)).
\end{equation}
Furthermore, in first order logic,  
\begin{equation}\label{eq:4} 
\forall r(\alpha\to\beta)\to (\exists r\alpha\to\exists r\beta)
\end{equation}
 is valid.
Let  then state the following
\begin{equation}\label{giacomo}
\phi(x,r)=  (x\in r\iff x\notin x )
\end{equation}
so that $\phi(x,r)$  is free from quantifiers and  (\ref{eq:2}) can be applied, so that
$$
\forall x(x\in r\iff x\notin x ) \to (r\in r\iff r\notin r).
$$
Let then also
$\alpha= \forall x\phi(x,r)$ e $\beta=\phi(r,r)$.  As we can notice from (\ref{eq:2}) we get 
$\forall r(\alpha\to\beta)$  and from (\ref{eq:1}) we draw  $\exists r\alpha$,  therefore from (\ref{eq:4}) applying two times modus ponens
we have $\exists r \beta$, that is 
\begin{equation}\label{eq:5} 
\exists r (r\in r\iff r\notin r).
\end{equation}

Since $r\in r\iff r\notin r$  is the negation of a propositional tautology
\begin{equation}\label{eq:6} 
(r\in r\iff r\notin r)\to\bot
\end{equation}
where $\bot$ refer to false, so by generalization
\begin{equation}\label{eq:7} 
\forall r((r\in r\iff r\notin r)\to\bot)
\end{equation}
or for the laws of quantifiers  
\begin{equation}\label{eq:8} 
\exists r(r\in r\iff r\notin r)\to\bot
\end{equation}
By modus ponens from (\ref{eq:5}) and (\ref{eq:8})  we can conclude
\begin{equation}\label{eq:9}
 \bot
\end{equation}
that is  (\ref{eq:1})  proves the false.
\end{proof}
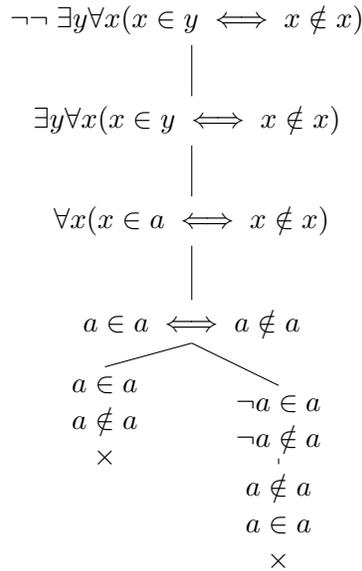
\begin{figure}[h]
\begin{center}
\begin{tikzpicture}[level distance=38pt, sibling distance=6em,
  every node/.style = { align=center}]
  \node {$ \neg \neg \; \exists y \forall x ( x \in y \iff x \notin x )$ }
child { node {$\exists y \forall x ( x \in y \iff x \notin x )$   }
child{ node {$ \forall x ( x \in a \iff x \notin x )$}
child{ node {$ a \in a \iff a \notin a$}
      child { node {$a\in a$\\$a\notin a$ \\ $\times$ }
       }
      child { node {$ \neg a \in a $\\ $ \neg a \notin a $}  
       child {node { $a \notin a $\\ $  a \in a $ \\$ \times$} }}}}};
\end{tikzpicture}
\end{center}
\caption[]{Semantic tree proving (\ref{antirussell}) }\label{albero1}
\end{figure}
We can show a  further proof of the antinomy  through the semantic tree method  as in
Figure (\ref{albero1}) 
(\cite{mendelson}, 141)(\cite{varzi}). 
By theorems  (\ref{RA}), (\ref{gl}), and the tree in Figure (\ref{albero1}),
definitely in first order logic  
 \begin{equation}\label{antirussell}
\vdash \; \neg \exists y \forall x ( x \in y \iff x \notin x ).
\end{equation}

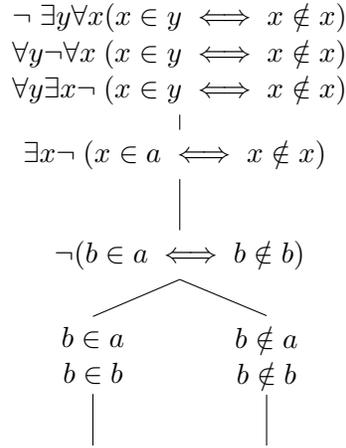
\begin{figure}[ht]
\begin{center}
\begin{tikzpicture}[level distance=38pt, sibling distance=6em,
  every node/.style = { align=center}]
  \node {$ \neg \; \exists y \forall x ( x \in y \iff x \notin x )$ 
  \\ $ \forall y \neg \forall x  \; ( x \in y \iff x \notin x ) $
  \\ $ \forall y \exists x  \neg \; ( x \in y \iff x \notin x ) $
  }
child { node {$ \exists x  \neg \; ( x \in a \iff x \notin x )$ }
child{ node {$ \neg ( b \in a \iff b \notin b)$}
      child { node {$b\in a$\\$b\in b$  }  child {node { $   \; $}}
       }
      child { node {$ b \notin a $\\ $  b \notin b  $}  
       child {node { $   \; $} }
       }}};
\end{tikzpicture}
\end{center}
\caption[]{A semantic tree for (\ref{russell}) }\label{albero2}
\end{figure}
This theorem is a negation of the existence of an object like the set of all those sets not belonging to themselves.  The proof of existence is  part of one of the most significant debates in contemporary logic, which is between proponents of solely existential proofs and others who believe that a demonstration of existence has to be constructive. In order to accept a proof that objects of a certain type exist, the latter require that a procedure be provided to produce them, while the former do not consider this requirement essential.  A classic example of a constructive theorem is the result already known to Euclid around 300 BC: given a prime number, there exists a larger one (which implies that prime numbers are infinite) (\cite{foundations}, 231).  An example of a purely existential demonstration can be found in the much more recent (19th century) Lemma of Bolzano Weierstrass: every infinite bounded set has at least one limit point (or accumulation point) (\cite{primer}, 45-48).  We are faced with two profoundly different conceptions of existence in mathematics.  The intuitionistic, or more generally constructivist, conception requires much more than the purely existential conception: one can admit that an object \emph{exists} only if a procedure is given that allows one to produce it step by step and in a finite number of steps. According to the classical logical conception, instead, to conclude that a mathematical object of a certain type \emph{exists} is sufficient to show that by introducing it we do not generate contradictions. As we can notice, perfectly aligned with this classical conception,  (\ref{antirussell})  is a theorem of nonexistence based directly on asserting a contradiction.   And this happens instead of concluding that $ \exists y \forall x ( x \in y \iff x \notin x )$ is an object that simply does not exist.  It is a mathematical object whose existence may make sense to assert as a hypothesis in a  proof  toward contradiction, the classical assumption in a reductio ad absurdum. But does it make sense to assert its existence in a theory when it has already been established that it does not exist?
If, indeed, we state in a  first order logic inferential system 
\begin{equation}\label{russell}\tag{R}
 \exists y \forall x ( x \in y \iff x \notin x ),
\end{equation}
from the standpoint of classical two-valued logic this is equivalent to assert a contradiction, a falsehood. 
Figure (\ref{albero2}) shows how the negation of statement (\ref{russell}) generates a semantic tree that has no closure, so for the semidecidability of the first-order logical calculus, it is not known whether a refutation is possible. 

Since it results in a contradiction of the kind $ \alpha \wedge \neg \alpha $, it appears that (\ref{russell}) is rejected by any logical system, thanks to the proof about (\ref{antirussell}).
We might therefore ask why we should persistently insist on affirming (\ref{russell}) within first-order inferential systems, inevitably making the system itself incoherent?
In its deductive complexity, if the system can generate (\ref{antirussell}) then it might also be able to demonstrate to us how to preserve its own consistency, that is, how, by asserting (\ref{russell}), we do not infer any contradiction or falsehood from it. Therefore, we might investigate if there is another inferential route that the system can take.
We also recall that theorem (\ref{antirussell}) does not only concern sets and the membership relation but also holds at first order in general, i.e. for every arbitrary binary relation $R$, $\; \vdash  \neg \exists y \forall x ( R(x, y) \iff \neg R(x, x) )$ (\cite{saeed}).
The following sections are also in this sense a deepening of the connection between logic and set theory.
We have recalled three distinct forms of what is currently accepted as evidence of Russell's antinomy since every distinct notation reveals novel insights that 
may shed light on these issues.

\bigskip

\section{Another inferential path}\label{sez1}
  In the ensuing discussion, we will revise the proof of theorem  (\ref{gl}) for a comprehensive inferential analysis before presenting a comparison with the other two kinds of proofs.
\begin{lemma}\label{nisba}
If   $\phi(x,r)=x\in r\iff x\notin x$  (\ref{giacomo})  then 
\begin{equation}\label{not_phi_rr}
 \forall x\: \phi (x,r) \to  \neg  \phi  (r,r).
 \end{equation}
 \end{lemma}
 \begin{proof} 
  $r\in r \iff r\notin r$  is the negation of a propositional tautology, then  $\neg (r\in r\iff r\notin r)$  is a tautology.  Given  (\ref{giacomo}) we have 
  $$ \forall x\: \phi (x,r) \to  \neg (r\in r \iff r\notin r)$$   so
 (\ref{not_phi_rr}) is always true.
 \end{proof} 
 A further proof dislplayed  in Figure (\ref{albero-nisba}).
 We also observe that the same kind of demonstration holds for $\gamma$ as  any antecedent in
$ \gamma \to  \neg  \phi  (r,r)$. If   $\phi(x,r)=x\in r\iff x\notin x$  (\ref{giacomo})  then 
\begin{equation}\label{neg_phi_rr}
 \gamma \to  \neg  \phi  (r,r).
 \end{equation}
 is always true.  
 Therefore, let's make a revision to theorem  (\ref{gl}), in accordance with Lemma (\ref{nisba}).
\begin{theorem}\label{lemma1}
Stating (\ref{russell}) may also not lead to contradiction (not even in classical logic).
\end{theorem}
\begin{proof} Let us suppose that (\ref{eq:1}) is true. 
In first order logic, once given (\ref{giacomo}), provided the substitution of $x$  by $r$ is valid, but it certainly is if  $\phi$ has no quantifiers, thanks to previous Lemma,  (\ref{not_phi_rr}) holds. 
Now for the generalization rule follows
\begin{equation}\label{eq:not_3} \forall r(\forall x\phi(x,r) \to \neg \phi(r,r)).
\end{equation}
We also have
\begin{equation}\label{eq:not_4} \forall r(\alpha\to\beta)\to (\exists r\alpha\to\exists r\beta).
\end{equation}
  Let us consider (\ref{giacomo}) so that $\phi(x,r)$ is quantifier free and we can apply (\ref{eq:not_3}). Let also
$\alpha=(\forall x\phi(x,r))$ and $\beta=\phi(r,r)$.
We can notice that  (\ref{eq:not_3}) yields  
$(\alpha\to \neg \beta)$. By generalization we get   $\forall r (\alpha\to \neg \beta)$. Affirming $ \exists r\ \forall x( x\in r\iff x\notin x)$  (\ref{eq:1}), i. e.  $\exists r \alpha$, 
and from 
$\forall r(\alpha\to\neg \beta)\to (\exists r\alpha\to\exists r\neg \beta)$ 
we obtain  $\exists r(\neg \beta)$,  that is  
 \begin{equation}\label{no_R}
 \neg(r \in r \iff r \notin r),
 \end{equation}
  which is the negation of false,  i.e. is true. That is, no contradiction.
\end{proof}
\begin{figure}[ht]
\begin{center}
\begin{tikzpicture}[level distance=38pt, sibling distance=6em,
  every node/.style = { align=center}]
  \node {$ \neg [\forall x \phi(x,r) \to  \neg ( r \in r \iff r \notin r )]$ }
child { node {$\forall x \phi(x,r) $ \\ $ \neg \neg ( r \in r \iff r \notin r )$ }
child{ node {$ ( r \in r \iff r \notin r )$}
      child { node {$r\in r$\\$r\notin r$ \\ $\times$ }}
      child { node {$ \neg r \in r $\\ $ \neg r \notin r $}  
       child {node { $r \notin r $\\ $  r \in r $ \\$ \times$} }}}};
\end{tikzpicture}
\end{center}
\caption[]{Semantic tree proving (\ref{nisba}) }\label{albero-nisba}
\end{figure}
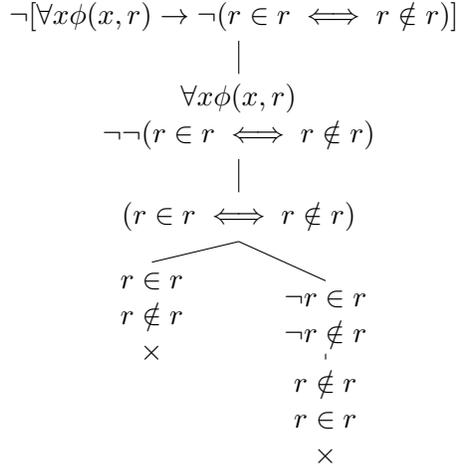
\begin{remark}
We can now clearly observe that (\ref{not_phi_rr}) in Lemma (\ref{nisba}) has the opposite consequent to that in (\ref{eq:2}), and that in the proof of Theorem (\ref{gl}), (\ref{giacomo}) directly conflicts with the application of the Predicative Axiom. The application of a substitution, namely  the substitution of $x$ by $r$, is also evident in definition (\ref{giacomo}).  Indeed, Lemma (\ref{nisba})  directly calls into question the applicability of the Predicative  Axiom to obtain (\ref{eq:2}).
\end{remark}
We also observe that since $\phi(r,r)$ is the negation of a propositional tautology, not only is (\ref{neg_phi_rr}) valid, but 
\begin{equation}\label{assert_phi_rr}
  \phi  (r,r) \to \gamma
 \end{equation}
  is also valid.
  $\gamma$ in (\ref{neg_phi_rr}) and (\ref{assert_phi_rr}), stands for any statement, i.e. for something  whose semantic content we do not know or do not discern. 
  
  Next section deals with first-order logic with identity, showing that there are laws that allow us to explore this unknown, finding a way to discern between the sets  $x$ and $y$ in (\ref{russell}), namely between the set of sets not belonging to themselves ad its members. We will deal with  the Law of the Indiscernibility of Identicals, $\forall x \forall y [ x = y \to (\phi (x) \iff \phi (y))]$,  
and how this law can help us learn to discern between sets.
For the next section, we will follow  Chang Keisler  using the symbol `$\equiv$' for the identity between terms (\cite{CK}).  For derivations from the Identity Laws the reader can refer to
(\cite{foundations}, 25), Ch. 8 in \cite{mendelson} and Ch. III in \cite{tarski}; see also Basic Law (III) in  (\cite{frege}, 20), and (\cite{cook}, A-26) for the two Laws, Indiscernibility of Identicals and Identity of Indiscernibles, as consequences of Basic Law (III).

\bigskip

\section{(In)-Discernibility}\label{sez2}
It is possible to demonstrate that logical systems that allow formalizing expressions  of type (\ref{russell}),  but prove (\ref{antirussell}),  such as theorems (\ref{RA}), (\ref{gl}), and the tree in figure (\ref{albero1}), can be considered as lacking something about  identity, as they neglect or ignore the Laws of Identity and  their implications.
 In other words, within these systems are allowed inferences that mystify    the validity of replacement in terms of the quantifiers and
 the validity of substitution in terms of the rules of identity. \bigskip

\begin{figure}[th]
\begin{center}
\begin{tikzpicture}[level distance=34pt, sibling distance=6em,
  every node/.style = { align=center}]
  \node {$ \neg \; [\forall x  ( \phi(x,r) \to \neg (x \equiv r)) \to \neg \phi(r,r)]$ 
  }
child { node {$ \forall x  ( \phi(x,r) \to \neg (x \equiv r)) $ }
child{ node {$\phi(r,r)$}
child{ node {$\phi(r,r) \to  \neg(r \equiv r) $}
      child { node {$\neg \phi(r,r)$ \\ $\times$ } 
       }
      child { node {$ \neg (r \equiv r) $ \\ $\times $}
      }}}};
\end{tikzpicture}
\end{center}
\caption[]{Semantic tree for (\ref{paola}) }\label{albero4}
\end{figure}
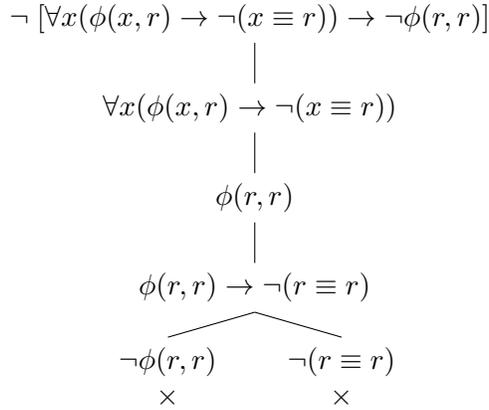

\begin{lemma}\label{LI}
Given   $         
\phi(x,r)=x\in r\iff x\notin x
$  (\ref{giacomo}) and  the  Laws of Identity  then $$ \forall x (\phi (x,r) \to  \neg (x \equiv r ) ) \to \neg \phi  (r,r). $$ 
\end{lemma}
\begin{proof} 
At the axiomatic deductive level, considering the Identity Axioms   (see 1.3.7 \cite{CK} 25), we have both $x \equiv y \rightarrow ( \phi (x,y)  \rightarrow \phi(x,x)  )$ and $x \equiv y \rightarrow ( \phi (x,x)  \rightarrow \phi(x,y)  )$, and thus also
\begin{equation}\label{ext1}
x \equiv r \rightarrow ( \phi (x,r)  \iff \phi (x,x)  )
\end{equation}
whence
 \begin{equation}
x \equiv r \rightarrow (x \in r  \iff x \in x  ).
\end{equation}
By the counterposed
 \begin{equation}
\neg (x \in r  \iff x \in x  )  \; \rightarrow \; \neg (x \equiv r ),
\end{equation}
 then
\begin{equation}\label{ext}
 (x \in r  \iff x \notin x  )  \rightarrow \neg (x \equiv r ).
\end{equation}
From   $\phi(x,r)=x\in r\iff x\notin x$    (\ref{giacomo}) and from (\ref{ext})  by modus ponens
\begin{equation}\label{noeq}
\neg (x \equiv r ).
\end{equation}
Which clearly means, since (\ref{giacomo}) defines $r$ as a set whose members do not belong to themselves, then, of course, thanks to the Laws 
of Identity, $r$ is not equal to $x$.  
 Whenever (\ref{giacomo}) then  always 
  $$\phi (x,r) \to  \neg (x \equiv r ),$$ and by the Generalization rule
\begin{equation}\label{p}
\forall x (\phi (x,r) \to  \neg (x \equiv r )).
\end{equation}
Thanks to (\ref{noeq}), $x$ and  $r$ are not identical,  $x$ cannot be replaced by  $r$ within (\ref{giacomo}), and consequently
\begin{equation}\label{paola}
\forall x (\phi (x,r) \to  \neg (x \equiv r )) \to \neg \phi  (r,r).
\end{equation}
\end{proof} 
\begin{figure}[ht]
\begin{center}
\begin{tikzpicture}[level distance=38pt, sibling distance=6em,
  every node/.style = { align=center}]
  \node {$ \neg  \; [ ( x \in y \iff x \notin x )  \to \neg (x \equiv y)]$ }
child { node {$( x \in y \iff x \notin x ) $   }
child{ node {$\neg \neg (x \equiv y)$\\$x \equiv y$}
child{ node {$ x \in x \iff x \notin x$}
      child { node {$x\in x$\\$x\notin x$ \\ $\times$ }
       }
      child { node {$ \neg x \in x $\\ $ \neg x \notin x $}  
       child {node { $x \notin x $\\ $  x \in x $ \\$ \times$} }}}}};
\end{tikzpicture}
\end{center}
\caption[]{Semantic tree proving  (\ref{ext}) }\label{albero3}
\end{figure}
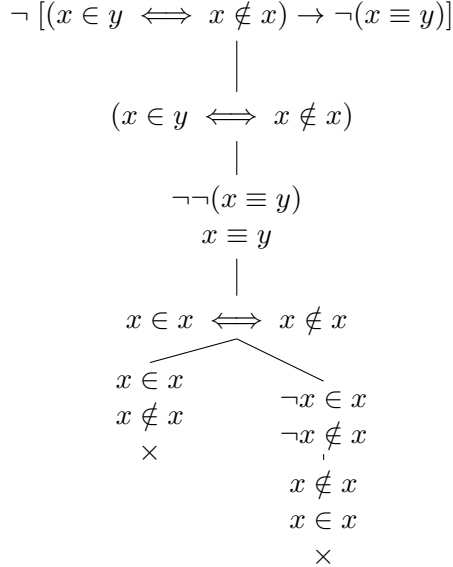
We can validate and at the same time explain Lemma (\ref{LI})  by showing the two semantic trees with applications of the identity rules  in figures (\ref{albero3}) and  (\ref{albero4}), respectively at nodes $ x \in x \iff x \notin x$ and $ \neg (r \equiv r) $.   For semantic trees with identities, see (\cite{varzi}).  About the tree in figure (\ref{albero3}), we notice that if we had used the constants $a$ and $b$  instead of respectively the variables $x$ and $y$, we would always have obtained a closed tree, namely the formula also applies independently of  the quantifiers. 
The tree also shows that (22) is true in   itself. 
That is, it holds both independently of and as a consequence of the axioms, as demonstrated by the deductions from (\ref{ext1}) to (\ref{ext}). 
Consequently, it generally holds in first-order logic 
\begin{equation}\label{senza_ext}
\vdash \forall x\forall y [ (x \in y \iff  x \notin x) \to  \neg (x \equiv y))].
\end{equation}
The semantic tree of figure (\ref{albero4}) displays further that
$$\vdash \forall x  ( \phi(x,r) \to \neg (x \equiv r)) \to \neg \phi(r,r),$$
which  confirms (\ref{paola}).

\medskip

\begin{remark}
As showed by Lemma  (\ref{LI}), when  (\ref{giacomo}), the premise of (\ref{paola}) is true, then $\neg \phi(r,r)$ is also true. 
With Lemma (\ref{LI}) and 
 figure (\ref{albero4}), the system offers an inferential investigation of  (\ref{not_phi_rr}) and Lemma (\ref{nisba}).
 \end{remark}
 
 \medskip
 
\begin{theorem}\label{magga}
Thanks to the Laws of Identity, in first order logic, no contradiction arises in affirming (\ref{russell}).
\end{theorem}
\begin{proof}
Suppose that (\ref{eq:1}) is true, then from (\ref{p}) and (\ref{paola}) it follows  $\neg  \phi(r,r)$, i.e. $ \neg(r \in r \iff r \notin r)$, i.e., the negation of the false, which is the true. No contradiction.
\end{proof}

\bigskip

\section{Uniqueness}
Having reached this point, we can ask ourselves, exactly what does the problem of Russell's antinomy consist of? Why did Frege write his famous appendix in Grundgesetze? 
As well known, the potential for expressing antinomy through Basic Law (V) appeared to be the cause of the system's incoherence (\cite{frege}).
Frege himself attempted an amendment to the Basic Law  (V), obtaining an inferential path, which did not lead to  contradiction in the expressions of the antinomy in Grundgesetze.
For a long time in the related literature the finger has been pointed at the axiom of Comprehension, which derives from Frege's Basic Law (V), from which (\ref{russell}) can be straight expressed.
According to theorems (\ref{lemma1}) and (\ref{magga}), with due consideration of Extensionality and Basic Law (III), even respectively in Zermelo-Fraenkel set theory   and Frege's Grundgesetze, given Russell's antinomy, the contradiction does not arise    (\cite{catta1}\cite{catta2}).

Let us recall in Zermelo-Fraenkel set theory the axiom of \emph{Comprehension }
$$\forall z_{1}...\forall z_{n}\exists y\forall 
x(x \in y \iff \varphi(x)),$$
where $ \varphi (x)$ is any formula in \emph{ZF},  $z_{1},...,z_{n}$ are the free variables of $\varphi (x)$ other than $x$, 
and $y$ is not a free variable of $\varphi(x)$; 
and the axiom of  \emph{Extensionality}
$$  \quad \forall x \forall y[  \forall z (z \in x \iff z \in y) \iff  
 x = y ],
$$
 (\cite{foundations} 27-30) (\cite{catta1} 3-6) (\cite{larisa} 109). 
 
 \begin{figure}[h]
\begin{center}
\begin{tikzpicture}[level distance=16mm, 
level 1/.style={sibling distance=100mm},
level 2/.style={sibling distance=70mm},
level 3/.style={sibling distance=40mm},
level 4/.style={sibling distance=40mm},
  every node/.style = { align=center}]
  \node {$ 
   \neg [\forall y \forall x [ ( x\in y \iff x \notin x )  \to \neg (x = y)]  \to    \neg \exists y \forall x (x \in y \iff x\notin x)    ]$
  }
child { node {$\forall y \forall x [ ( x\in y \iff x \notin x )  \to \neg (x = y)]  $ \\ $  \neg\neg \exists y \forall x (x \in y \iff x\notin x)   $ }
      child { node {$\exists y \forall x (x \in y \iff x\notin x)$\\ $ \forall x (x \in a \iff x\notin x)$\\$(a \in a \iff a\notin a)$ }
      child {node {$( a\in a \iff a \notin a )  \to \neg (a = a)$}
      child{ node {$\neg ( a\in a \iff a \notin a ) $\\ $\times$}} child{node {$\neg (a = a) $\\ $\times$}} }  }
      };
\end{tikzpicture}
\end{center}
\caption[]{Semantic tree expanding the tree in Figure (\ref{albero4}) }\label{albero7}
\end{figure}
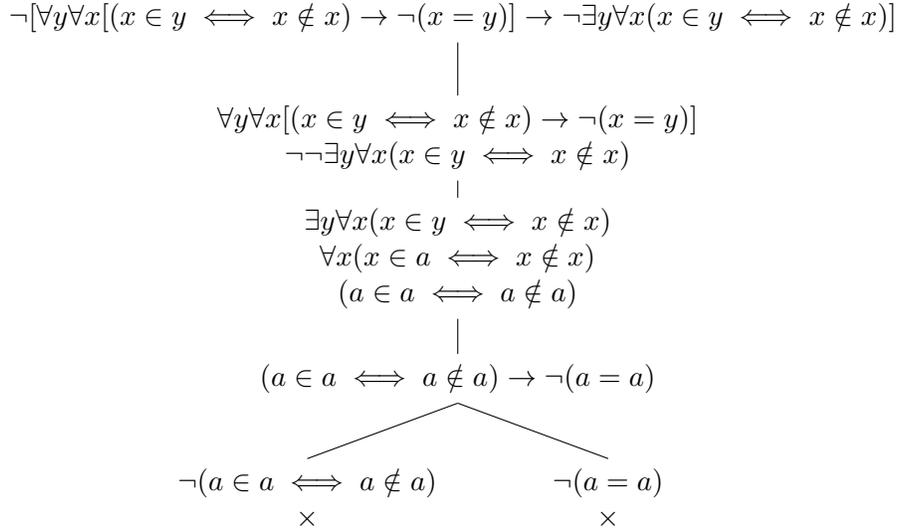

 \medskip
 
When we assert (\ref{russell})  in \emph{ZF}, we are once again dealing with $\gamma$ as in (13) and (17), we can get to $\gamma  \to \neg$(\ref{russell})  and (\ref{russell}) $\to \gamma$, but in \emph{ZF} we are allowed to  explore $\gamma$ further thanks to the Theory of Definition.

The tree that does not close in Figure  (\ref{albero2}) supports the belief, dating back to the well-known appendix of Frege's Grundgesetze, that the formal assertability of  (R) entails the inconsistency of the theory or system that enables it, which implies that we do not know how to deal with (\ref{russell}). 
Figure  (\ref{albero7}) illustrates how, taking Extensionality into account,
we  know how to block the inference of contradiction when we state  
 (\ref{russell}), and thus how not to damage the coherence of the system.
The tree in Figure (\ref{albero7}) could be closed after the node  ``$a \in a \iff a \notin a$", similarly to the tree in  Figure (\ref{albero1}), but we chose to extend it to display the second possible semantic path associated  with Extensionality.  The tree  in Figure  (\ref{albero7}) can be regarded as a \emph{ZF} expansion of the tree in Figure  (\ref{albero4})  further shedding light on the contradictory nature of $\phi (r,r)$ and its inferential connections.

 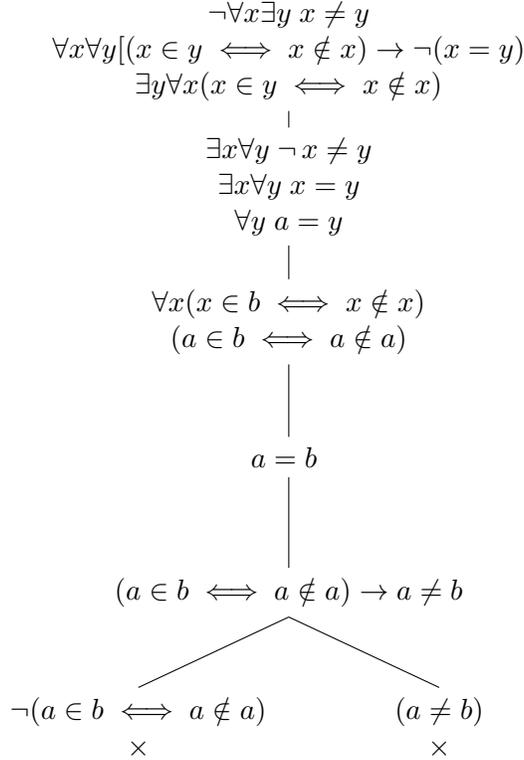
\begin{figure}[h]
\begin{center}
\begin{tikzpicture}[level distance=18mm, 
level 1/.style={sibling distance=100mm},
level 2/.style={sibling distance=70mm},
level 3/.style={sibling distance=40mm},
level 4/.style={sibling distance=40mm},
  every node/.style = { align=center}]
  \node {$ \neg \forall x \exists y \; x \neq y$ \\
  $ \forall x \forall y [ ( x\in y \iff x \notin x )  \to \neg (x = y) $ \\  
   $ \exists y \forall x (x \in y \iff x\notin x)    $
  }
child { node {$  \exists x \forall y \; \neg \, x \neq y$ \\ $ \exists x \forall y \;  x = y  $ \\ $\forall y \;  a = y$ }
      child { node {$ \forall x (x \in b \iff x\notin x)$\\ $( a\in b \iff a \notin a )$ }
      child { node {$ a = b$ }
      child {node {$( a\in b \iff a \notin a ) \to a \neq b$}
      child{ node {$\neg ( a\in b \iff a \notin a ) $\\ $\times$}} child{node {$ (a \neq b) $\\ $\times$} } }  }
      }  };
\end{tikzpicture}
\end{center}
\caption[]{Semantic tree for Theorem (\ref{maggaZF})}\label{urka}
\end{figure}

The Axiom of Comprehension can be considered as a schema originally conceived for  defining sets.
Because of its role  in introducing new sets, i.e. in  bringing into existence  ever new mathematical objects,   this certainly concerns the rules for introducing new symbols into already well-established theories. 
We refer here to the rules of the Theory of Definition based on the two Criterions of  Eliminability  and Non-Creativity, originally by S. Le\'{s}niewsky (\cite{catta1}, \cite{suppes}).  In particular, by the rule for defining a \emph{new operation symbol} (or a new individual 
constant, i.e., an operation symbol of rank zero), we are required to have a preceding theorem 
which guarantees that the  operation is uniquely defined.
If the 
restriction on the uniqueness is dropped then a contradiction can be 
derived  (\cite{suppes}, 159) (\cite{catta1}, \cite{catta2}).

Extensionality, without any additional axioms,  implies that for 
every condition $\varphi(x)$ on $x$ (in Comprehension) there exists ``at 
most" one 
set $y$ which contains \emph{exactly} those elements $x$ which fulfil the 
condition 
$\varphi(x)$ (\cite{foundations},  31), (\cite{catta1}).
In other words, if  there is a set $y$ such that $\forall 
x(x \in y \leftrightarrow \varphi(x))$, $y$ is unique. It can be shown as 
follows. If $y'$ is also such, i.e. $\forall 
x(x \in y' \leftrightarrow \varphi (x))$, then we have, obviously $\forall 
x(x \in y' \leftrightarrow x \in y)$, and then by Extensionality, $y' = y$. 
(\cite{foundations} 31). 
Noticeably, the proof in Theorem  (\ref{RA}) involve only the axiom of Comprehension.  As is easily demonstrated, the contradiction in its proof does not occur when  defining (\ref{russell}) if Extensionality is used in addition to Comprehension (\cite{catta1}, \cite{catta2}, \cite{slai}). In set theory, Extensionality is a derivation of the Laws of Identity, which in turn are  in Frege's Grundgesetze a derivation from Basic Law (III) (\cite{catta1} \cite{catta2}).  By Extensionality the 
above proof of Theorem  (\ref{RA}) can not be concluded. 
\begin{theorem}\label{maggaZF}
Thanks to the Axiom of Extensionality, in Zermelo-Frankel set theory, no contradiction arises in affirming (\ref{russell}).
\end{theorem}
\begin{proof}
If, as an example of   
Comprehension, we define
\begin{equation}\label{comp1}
x \in y \iff x \notin x
\end{equation}
then by Extensionality  we always obtain (Figure (\ref{albero3}) and (\cite{catta1} 3-6)),
\begin{equation}\label{star}
(x \in y \iff x \notin x) \to \neg(x = y),
\end{equation}
so by modus ponens
\begin{equation}\label{twostars}
x \neq y.
\end{equation}
Given (\ref{comp1}) and (\ref{star}), by  (\ref{twostars}) in the domain of interpretation $x$ and $y$ are not the same individual, so $y$ cannot be considered identical to $x$. Precisely, Figure (\ref{urka}), $  \forall x \exists y \; x \neq y$, so that for any $x$  there exists at least one $y$ such that  $ x \neq y$.
This is an applicative case of the uniqueness 
condition.
Indeed a special case because it involves negation and therefore 
complementation (wrt. (\ref{senza_ext})).
For this reason, the inference of the proof of Theorem  (\ref{RA})
cannot be terminated since $``x = y"$ cannot be assumed and we do not obtain $``y \in y 
\text{ if and only if }  y \notin y"$.
 (\cite{catta1}, \cite{catta2}, \cite{slai}).
 \end{proof}
 
As we can notice, comparing the arguments of the theorems (\ref{RA}) and  (\ref{gl}), (\ref{giacomo}) is an example of Comprehension formula like (\ref{comp1}), (\ref{star}) is  (\ref{ext}) and (\ref{twostars}) is   (\ref{noeq}).
Further, (\ref{paola}) is quite close to the consequences of Frege's amendment of  Basic Law (V)  (\cite{frege}), so to block in a quite similar way the derivation of the contradiction, with the only difference that  (\ref{paola}) is simply derived from the Laws of Indentity, i.e., from Basic Law (III), leaving Frege's Basic Law (V) untouched
 (for any details see \cite{catta1}, \cite{catta2}).

Frege amended Basic Law (V),  embedding the negation of the identity between the function defining the class of classes not belonging to themselves and its extension.
On a strictly deductive level, such a negation blocks the derivation of the contradiction, although as shown by S. Le\'{s}niewski, the amendment of the axiom is not formally satisfactory.   

As already demonstrated, there is no need of this restriction of the Basic  Law (V) ([2], [3]).  The restrictive clause instead of being included in the Basic Law (V), can simply be deduced from the Basic Law (III), in a manner quite analogous to what is shown in Lemma (\ref{LI}) and in Theorem (\ref{magga}). This prevents the contradiction from being derived, with inferences strikingly similar to those of Frege, but this happens by deducing the restriction from the axioms of the system, without modifying them.

\section{Ending Note}
There is a specific reason why, when introducing a new symbol into a theory, one must take into account uniqueness and therefore, as we have seen, implicitly the Laws of Identity.  In a first-order logic system, identity is  \emph{relative} to the language of the system, that is, identity is a relation that refers only to previously defined predicates. The uniqueness violations  implicitly mean opening the universe to objects with characteristics that are not yet well defined. In the case of infinite sets, each time we add a new element, we must require that the relations between the previous elements remain unchanged, that is, that they are not modified by the new addition.  As seen thus far, in a first-order set theory, it is essential to consider the uniqueness of each new set through Extensionality when Comprehension introduces some set defined by self-reference as in (\ref{russell}).

\end{document}